\documentclass[a4paper,12pt,reqno]{article}

\usepackage{amsthm}
\usepackage{amsmath,amssymb,latexsym,amsfonts,mathrsfs}
\usepackage[dvipdfmx]{graphicx}
\usepackage{bm}

\usepackage{color}
\usepackage{comment}
\usepackage{mathtools}
\usepackage{fancyhdr}
\usepackage[top=30truemm,bottom=30truemm,left=20truemm,right=20truemm]{geometry}
\usepackage{enumerate}

\newcommand{\R}{\mathbb{R}}
\newcommand{\C}{\mathbb{C}}
\newcommand{\N}{\mathbb{N}}
\newcommand{\Z}{\mathbb{Z}}

\newcommand{\SCR}[1]{{\mathscr #1}}

\newcommand{\CAL}[1]{{\cal #1}}

\newcommand{\J}[1]{\left\langle #1 \right\rangle}

\theoremstyle{plain}
\newtheorem{Thm}{{\bf Theorem}}[section]

\newtheorem{Lem}[Thm]{{\bf Lemma}}
\newtheorem{Prop}[Thm]{{\bf Proposition}}

\newtheorem{Ass}[Thm]{{\bf Assumption}}

\theoremstyle{definition}
\newtheorem{Def}[Thm]{{\bf Definition}}
\newtheorem{Rem}[Thm]{{\bf Remark}}

\newcounter{Exami}

\allowdisplaybreaks[2]

\makeatletter
 \def\address#1#2{\begingroup
 \noindent\parbox[t]{7.8cm}{%
 \small{\scshape\ignorespaces#1}\par\vskip0ex
 \noindent\small{\itshape E-mail}%
 \/: #2\par\vskip4ex}\hfill%
 \endgroup}%
\makeatother

\begin{document}
\fontencoding{T1}\selectfont
\title{Strichartz estimates in Wiener amalgam spaces for Schr\"{o}dinger equations with at most quadratic potentials}
\author{Shun Takizawa}
\date{\today}
\maketitle
\begin{abstract}
For Schr\"{o}dinger equations with potentials which grow at most quadratically at spatial infinity, we prove Strichartz estimates in Wiener amalgam spaces. These estimates provide a stronger recovery of  local-in-space regularity than the classical Strichartz estimates in Lebesgue spaces. Our result is a generalization of the results on Strichartz estimates in Wiener amalgam spaces by Cordero and Nicola, which are stated for the potentials $V(x)=0, \pm |x|^2/2$.
\end{abstract}
\section{Introduction} 
In this paper, we study the following Cauchy problem for Schr\"{o}dinger equations with potentials:
\begin{align}\label{CP1}
\begin{cases}
&i\partial_{t}u(t,x)= -\frac{1}{2}\Delta u(t,x)+V(x)u(t,x), \hspace{3mm}(t,x)\in \R \times \R^n,
\\&u|_{t=0}=u_{0}\in L^2(\R^n).
\end{cases}
\end{align}
where $u(t), u_0$ are $\C$-valued and $V$ is $\R$-valued.
We impose the assumption that the potential $V$ has at most quadratic growth.
\begin{Ass}\label{AssV}
The potential is $V\in C^{\infty}(\R^n; \R)$ and there exists $C_{\alpha}>0$ such that
\begin{equation*}
|\partial_{x}^{\alpha}V(x)|\leq C_{\alpha},\hspace{3mm}x\in \R^n, \hspace{2mm}|\alpha|\geq2.
\end{equation*}
\end{Ass}
We put $H=-\frac{1}{2}\Delta+V$. Under Assumption \ref{AssV}, it is known that the operator $H$ generates the propagator $U(t)=e^{-itH}$, which is a one parameter unitary group and the solution to (\ref{CP1}) is given by $u(t)=e^{-itH}u_0$ (see e.g. \cite{Faris-Lavine, Kato1}). 

Let us recall the global or local-in-time (homogeneous) Strichartz estimates for the Schr\"{o}dinger equation which state that 
\begin{align} \label{Stri1}
\|e^{-itH}u_0\|_{L^r(\R; L^p(\R^n))}\leq C \|u_0\|_{L^2(\R^n)} 
\end{align}
or for all $T>0$,
\begin{align} \label{Stri2}
\|e^{-itH}u_0\|_{L^r([-T,T]; L^p(\R^n))}\leq C_T \|u_0\|_{L^2(\R^n)},
\end{align}
where $(p,r)$ fulfills
\begin{align*}
2\leq p,r \leq \infty, \hspace{3mm}\frac{n}{p}+\frac{2}{r}=\frac{n}{2},\hspace{3mm} (n,p,r)\ne (2,\infty,2).
\end{align*}
These estimates express a gain of local regularity in $x$ of the solution $u(t)$ for almost all $t$ since $p\geq2$. 
For the case $V=0$, the global-in-time Strichartz estimate (\ref{Stri1}) has been established by Strichartz \cite{Strichartz} (Case: $r=p$), Ginibre-Velo \cite{Ginibre-Velo1}, Yajima \cite{Yajima4} (Case: $r>2$) and Keel-Tao \cite{Keel-Tao} (Case: $r=2$). Subsequently the Strichartz estimates have been extended to the case with potentials which decay at spatial infinity (see \cite{Rod-Sch, BPST1, BPST2, Goldberg, Mizutani, Taira} and the references therein).
On the other hand we remark on the case with potentials which have at most quadratic growth at spatial infinity. Fujiwara \cite{Fujiwara2} has shown that under Assumption \ref{AssV} the fundamental solution has $|E(t,x,y)|\leq C |t|^{-n/2}$ for all $x,y\in \R^n$ and $t\ne0$ sufficiently small,  
where the fundamental solution is defined by the integral kernel of solution $u(t)$ to (\ref{CP1}):  $u(t,x)=\int_{\R^n} E(t,x,y) u_0(y) dy$. The estimate of $E(t,x,y)$ leads short-time dispersive estimate:
 $\|u(t)\|_{L^\infty}\leq C |t|^{-n/2} \|u_0\|_{L^1}$, which implies the  local-in-time Strichartz estimate (\ref{Stri2}) by the $TT^{*}$-method (see \cite{Ginibre-Velo1, Keel-Tao}).   
 Next we recall the (homogeneous) Strichartz estimates in Wiener amalgam spaces for (\ref{CP1}) which are established by Cordero-Nicola \cite{Cordero-Nicola1, Cordero-Nicola2} and state that
\begin{align} 
&\|e^{it\Delta/2}u_0\|_{W(L^{r/2}, L^r)_{t} W(\SCR{F}L^{p'}, L^p)_{x}} \leq C \|u_0\|_{L^2(\R^n)},  \label{Stri3}
\\
&\|e^{it(\Delta+|x|^2)/2}u_0\|_{W(L^{r/2}, L^r)_{t} W(\SCR{F}L^{p'}, L^p)_{x}} \leq C \|u_0\|_{L^2(\R^n)}  \label{Stri4}
\end{align}
and for all $T>0$,
\begin{align} \label{Stri5}
\|e^{it(\Delta-|x|^2)/2}u_0\|_{L^{r/2}([-T,T];W(\SCR{F}L^{p'}, L^p)) }\leq C_T \|u_0\|_{L^2(\R^n)}
\end{align}
with the $(p,r)$ satisfying the admissible condition:
\begin{align*}
2\leq p \leq \infty,\hspace{3mm} 4< r \leq \infty,\hspace{3mm}\frac{n}{p}+\frac{2}{r}=\frac{n}{2}.
\end{align*}
Here, $\|f(t,x)\|_{W(A,B)_{t} W(C,D)_{x}}:=\| \|f(t,x)\|_{W(C,D)_{x}} \|_{W(A,B)_{t}}$ and $W(A,B)$ is Wiener amalgam spaces, which are defined by Definition \ref{DefWA} below. The estimates (\ref{Stri3})-(\ref{Stri5}) imply that the solution $u(t)=e^{-itH}u_0$ belongs to $W(\SCR{F}L^{p'}, L^p)$ for almost all $t$. In addition Cordero-Nicola proved that (\ref{Stri3})-(\ref{Stri5}) also hold for the endpoint case : $(p,r)=(\frac{2n}{n-1}, 4), \hspace{2mm}n>1$ if we replace $\SCR{F}L^{p'}$ with $\SCR{F}L^{p',2}$, where $L^{p',2}$ is the one of the Lorentz spaces.
 Noting that $W(\SCR{F}L^{p'}, L^p)\hookrightarrow W(\SCR{F}L^{p',2}, L^p) \hookrightarrow L^p$ (see below Section 2.4), we find that the estimates (\ref{Stri3}), (\ref{Stri4}) and (\ref{Stri5}) yield stronger recovery of regularity in $x$ than the estimates (\ref{Stri1}) and (\ref{Stri2}), respectively. 
%
%
%

We give the terminology in order to state our result.
\begin{Def}
We call $(p,r)$ an admissible pair if $(p,r)$ satisfies 
\begin{align*}
2\leq p \leq \infty,\hspace{3mm} 4\leq r \leq \infty,\hspace{3mm}\frac{n}{p}+\frac{2}{r}=\frac{n}{2},\hspace{3mm} (n,p,r)\ne (1,\infty,4).
\end{align*}
The admissible pair $(p,r)=(\frac{2n}{n-1}, 4)$ is said to be endpoint.
\end{Def}
The following theorem is our main result and an extension of the results by Cordero-Nicola \cite{Cordero-Nicola1, Cordero-Nicola2} to the case with general potentials of at most quadratic growth. 
\begin{Thm}\label{thm1}
Suppose that $V$ satisfies Assumption \ref{AssV} and $T>0$. Let $(p,r)$ and $(\tilde{p}, \tilde{r})$ be non-endpoint admissible pairs. Then we have the homogeneous Strichartz estimates
\begin{align} \label{MainInq1}
\| e^{-itH}u_0\|_{L^{r/2}\left([-T,T]; W(\SCR{F}L^{p'}, L^p)\right)} \leq C_{T} \|u_0\|_{L^2},
\end{align}
the dual homogeneous Strichartz estimates
\begin{align} \label{MainInq2}
\left\| \int_{-T}^{T} e^{isH} F(s)ds \right\|_{L^2} \leq C_{T} \|F\|_{L^{(r/2)'}\left([-T,T]; W(\SCR{F}L^{p}, L^{p'})\right)},
\end{align}
and the retarded (or also inhomogeneous) Strichartz estimates
\begin{align} \label{MainInq3}
\left\| \int_{0}^{t} e^{-i(t-s)H}F(s) ds \right\|_{L^{r/2}\left([-T,T]; W(\SCR{F}L^{p'}, L^p)\right)} \leq C_{T} \|F\|_{L^{(\tilde{r}/2)'}\left([-T,T]; W(\SCR{F}L^{\tilde{p}}, L^{\tilde{p}'})\right)}.
\end{align}
Moreover for the endpoint case: $(p,r)=(\frac{2n}{n-1}, 4)$,$\hspace{2mm} n>1$, we obtain 
\begin{align} 
&\| e^{-itH}u_0\|_{L^{2}\left([-T,T]; W(\SCR{F}L^{p',2}, L^p)\right)} \leq C_{T} \|u_0\|_{L^2}, \label{MainInq4}
\\
&\left\| \int_{-T}^{T} e^{isH} F(s) ds\right\|_{L^2} \leq C_{T} \|F\|_{L^{2}\left([-T,T]; W(\SCR{F}L^{p,2}, L^{p'})\right)}. \label{MainInq5}
\end{align}
The retarded estimates (\ref{MainInq3}) still hold with non-endpoint $(p,r)$ and endpoint $(\tilde{p}, \tilde{r})$, if
 we replace $\SCR{F}^{\tilde{p}'}$ by $\SCR{F}^{\tilde{p}',2}$.  Similarly (\ref{MainInq3})  holds for endpoint $(p,r)$ and non-endpoint $(\tilde{p}, \tilde{r})$ if we replace $\SCR{F}L^{p'}$ by $\SCR{F}L^{p',2}$, Finally (\ref{MainInq3}) holds for endpoint $(p,r)=(\tilde{p}, \tilde{r})$ if we replace $\SCR{F}L^{p'}$ by $\SCR{F}L^{p',2}$ and  $\SCR{F}^{\tilde{p'}}$ by $\SCR{F}^{\tilde{p}',2}$, respectively.
\end{Thm}
\begin{Rem} \label{Rem}
Local-in-time Strichartz estimates of Theorem \ref{thm1} do not 
in general extend to global-in-time, as the case with harmonic oscillator : $V=\frac{1}{2}|x|^2$ explicitly shows. 
However in the case where the potential is the so-called Stark potential $V(x)=E\cdot x$ with given constant $E\in\R^n$, we can prove the global-in-time Strichartz estimates in Wiener amalgam spaces. Namely we have
\begin{align}
&\| e^{it(\Delta/2-E\cdot x)}u_0\|_{W(L^{r/2}, L^r)_{t} W(\SCR{F}L^{p'}, L^p)_{x}} \leq C \|u_0\|_{L^2}, \notag
\\
&\left\| \int_{-\infty}^{\infty} e^{-is(\Delta/2-E\cdot x)} F(s)ds \right\|_{L^2} \leq C \|F\|_{W(L^{(r/2)'}, L^{r'})_{t} W(\SCR{F}L^{p}, L^{p'})_{x}},  \notag
\\
&\left\| \int_{0}^{t} e^{-i(t-s)(\Delta/2-E\cdot x)}F(s) ds \right\|_{W(L^{r/2}, L^r)_{t} W(\SCR{F}L^{p'}, L^p)_{x}} \leq C \|F\|_{W(L^{(\tilde{r}/2)'}, L^{\tilde{r}'})_{t} W(\SCR{F}L^{\tilde{p}}, L^{\tilde{p}'})_{x}} \label{Stark}
\end{align}
for any non-endpoint admissible pairs $(p,r), (\tilde{p}, \tilde{r})$.
Moreover for the endpoint case: $(p,r)=(\frac{2n}{n-1}, 4)$,$\hspace{2mm} n>1$, we obtain 
\begin{align*} 
&\| e^{it(\Delta/2-E\cdot x)}u_0\|_{W(L^2, L^4)_{t} W(\SCR{F}L^{p',2}, L^p)_{x}} \leq C \|u_0\|_{L^2},
\\
&\left\| \int_{-\infty}^{\infty} e^{-is(\Delta/2-E\cdot x)} F(s)ds \right\|_{L^2} \leq C \|F\|_{W(L^2, L^{4/3})_{t} W(\SCR{F}L^{p,2}, L^{p'})_{x}},
\end{align*}
The retarded estimates (\ref{Stark}) still hold with non-endpoint $(p,r)$ and endpoint $(\tilde{p}, \tilde{r})$, if
 we replace $\SCR{F}^{\tilde{p}'}$ by $\SCR{F}^{\tilde{p}',2}$.  Similarly (\ref{Stark})  holds for endpoint $(p,r)$ and non-endpoint $(\tilde{p}, \tilde{r})$ if we replace $\SCR{F}L^{p'}$ by $\SCR{F}L^{p',2}$, Finally (\ref{Stark}) holds for endpoint $(p,r)=(\tilde{p}, \tilde{r})$ if we replace $\SCR{F}L^{p'}$ by $\SCR{F}L^{p',2}$ and  $\SCR{F}^{\tilde{p'}}$ by $\SCR{F}^{\tilde{p}',2}$, respectively.

This statement can be proved easily by combining the explicit formula of the fundamental solution and the results for the case $E=0$ by Cordero-Nicola \cite{Cordero-Nicola2}.
The proof is given in the Appendix.
\end{Rem}
\begin{Rem}
As an application of Theorem \ref{thm1}, we can study the well-posedness in $L^2$ for Schr\"{o}dinger equations with $H=-\frac{1}{2}\Delta+V$ as the free Hamiltonian: 
\begin{align} \label{CP2}
\begin{cases}
&i\partial_{t}u(t,x)+Hu= \widetilde{V}(t,x)u(t,x), \hspace{3mm}(t,x)\in \R \times \R^n,
\\&u|_{t=0}=u_{0}\in L^2(\R^n).
\end{cases}
\end{align}
Here, we suppose that the non-perturbed potential $V$ satisfies Assumption \ref{AssV} and the  perturbed potential $\widetilde{V}\in L^{r}([-T,T]; W^{p,p'}(\R^n))$ with $\frac{1}{r}+\frac{n}{p}\leq1, 1\leq r<\infty, n<p\leq\infty.$
Then, the Cauchy problem (\ref{CP2}) has a unique solution 
\begin{align*}
&u\in C([-T,T];L^2(\R))\cap L^{r/2}([-T,T];W^{p,p'}(\R))\hspace{3mm} \text{if} \hspace{3mm} n=1,
\\
&u\in C([-T,T];L^2(\R^n))\cap L^{r/2}([-T,T];W^{p,p'}(\R^n))\cap L^{2}([-T,T];W(\SCR{F}L^{\frac{2n}{n+1},2}, L^{\frac{2n}{n-1}})(\R^n)) 
\end{align*}
if $n\geq 2$ by the method of D'Ancona-Pierfelice-Visciglia \cite{DAncona}. This proof is the same argument as   
\cite[Theorem 6.1]{Cordero-Nicola3} (see also \cite{Cordero-Nicola1, Cordero-Nicola2}).
\end{Rem}
\subsection*{Idea and outline of the proof.}
The standard argument to prove Strichartz estimates is to deduce the dispersive estimates from the explicit formula or properties of the fundamental solution.
Indeed, if this is obtained, then together with the $L^2$ conservation law we can prove Strichartz estimates by the $TT^{*}$ method. Cordero-Nicola \cite{Cordero-Nicola1, Cordero-Nicola2} established the dispersive estimates of the Wiener amalgam space type:
\begin{align*} 
&\|e^{it\Delta/2}u_0\|_{W(\SCR{F}L^{1}, L^{\infty})} \leq C_n (|t|^{-n}+|t|^{-n/2})\|u_0\|_{W(\SCR{F}L^{\infty}, L^1)} \hspace{3mm} \text{for}\hspace{3mm} t\ne0, 
\\
&\|e^{it(\Delta-|x|^2)/2}u_0\|_{W(\SCR{F}L^{1}, L^{\infty})} \leq C_n |t|^{-n}\|u_0\|_{W(\SCR{F}L^{\infty}, L^1)}  \hspace{3mm} \text{for}\hspace{3mm} 0<|t|<\pi, 
\\
&\|e^{it(\Delta+|x|^2)/2}u_0\|_{W(\SCR{F}L^{1}, L^{\infty})} \leq C_n (|\sinh t|^{-n}+|\sinh t|^{-n/2}) \|u_0\|_{W(\SCR{F}L^{\infty}, L^1)} \hspace{3mm} \text{for}\hspace{3mm} t\ne0 
\end{align*}
by using the metaplectic representations, which can be applied only to $V(x)=k|x|^2, \hspace{3mm}k\in\R$. Moreover they proved (\ref{Stri3})-(\ref{Stri5}) by Keel-Tao's $TT^{*}$ method \cite{Keel-Tao} similarly to the Lebesgue space case from density and duality (see below Proposition \ref{Prop4}) but the method requires some additional technical ingredients, due to the different nature of the Wiener amalgam spaces. 

However, we cannot apply the method to the general potential $V$ fulfilling  Assumption \ref{AssV}.
In order to overcome this difficulty, we use the classical Hamiltonian flow combined with short-time Fourier transform techniques instead of the metaplectic representation.
In more detail, we transform (\ref{CP1}) into a Duhamel type formula $u(t)=U_0(t)u_0-i\int_{0}^{t}R(t,s)u(s)ds$, where $U_0(t)$ and $R(t,s)$ are corresponding to a parametrix operator and remainder, respectively (see also Section 5.1 in \cite{Takizawa}). Moreover, we attain the dispersive type estimates of the operators $U_0(t), R(t,s)$ by means of its factorizations as compositions of a change-of-variables operator along the Hamiltonian flow and the short-time Fourier transform.
An important feature of our proof is that it does not require the dispersive estimate in Wiener amalgam spaces for the original solution to (\ref{CP1}). 

 The rest of the paper is organized as follows. In Section 2, we introduce the notation and tools required for the proof, and prepare several lemmas concerning classical trajectories.
In Section 3, we present two key lemmas, from which the dispersive type estimates for $U_0(t)$ and $R(t,s)$  follow.
In Section 4.1, we transform (\ref{CP1}) into an integral equation: $u(t)=U_0(t)u_0-i\int_{0}^{t}R(t,s)u(s)ds$, which also appears in \cite{Takizawa} by making use of the method via Kato-Kobayashi-Ito \cite{KKI2}.
In Section 4.2 and 4.3, we derive the Strichartz type estimate for $U_0(t)$ and $R(t,s)$ by $TT^{*}$-method, respectively.
In Section 4.4, we prove (\ref{MainInq1}), (\ref{MainInq2}), (\ref{MainInq4}) and (\ref{MainInq5}).
In Section 4.5, we show (\ref{MainInq3}) and its endpoint case.
\section{Preliminaries}
\subsection{Notation}
Let $\nabla=(\partial_{x_1},\ldots,\partial_{x_n})$. $\J{x}$ stands for $(1+|x|^2)^{1/2}$.
We use the notation $\SCR{F}f(\xi)=\widehat{f}(\xi)=\int_{\R^n}f(x)e^{-ix\cdot\xi}dx$ for Fourier transform of $f$ 
and $\SCR{F}^{-1}f(x)=\check{f}(x)=(2\pi)^{-n} \int_{\R^n}f(\xi)e^{ix\cdot\xi}d\xi$ for the  inverse Fourier transform of $f$. We often write $\int$ instead of $\int_{\R^m}$ with $m\in \N$ for short.
We also often write $x^2:=|x|^2$ and $xy:=x\cdot y$ for $x, y\in \R^n$.
We denote $\CAL{S}$ by the Schwartz space and $\CAL{S}'$ by the set of all the tempered 
distributions. We denote $C_{b}^{\infty}(\R^n)$ by the set of all the functions of $C^{\infty}$-class and its all derivatives are bounded. 
We define the Fourier Lebesgue space by $\SCR{F}L^p=\{f\in \CAL{S}': \|f\|_{\SCR{F}L^p}:=\|\widehat{f}\|_{L^p}<\infty\}$.
For a Banach space $X$ and $X$-valued measurable functions $f$, we write the norm $\|f\|_{L^{p}X}=\left(\int \|f(x)\|_{X}^{p} dx \right)^{1/p}$. For Banach spaces $X$ and $Y$, we define $B(X, Y)$ the Banach space of bounded linear operators from $X$ to $Y$ and write the operator norm $\|\cdot\|_{X\to Y}:=\|\cdot\|_{B(X,Y)}$. The index $p'$ stands for the H\"{o}lder conjugate of $p\in[1,\infty]$, that is, $p'=p/(p-1)$.
Throughout this paper the letter $C$ denotes a constant, which may be different in each occasion.
Furthermore we often do not write subscript $n$ and $V$ when constants $C$ depend on the space dimension $n$ and the potential $V$.
We define the free Schr\"{o}dinger propagator $e^{it\Delta/2}$ by
\begin{align*}
e^{it\Delta/2}f(x)=\SCR{F}^{-1} e^{-it|\xi|^2/2} \SCR{F}f(x)=(2\pi it)^{-n/2}\int_{\R^n} e^{i\frac{(x-y)^2}{2t}} f(y) dy.
\end{align*} 
\subsection{Short-time Fourier transform (STFT)} 
We introduce the short-time Fourier transform used in the definition of the Wiener amalgam space.
\begin{Def}[STFT]
Let $g \in \CAL{S}(\R^n)\setminus \{0\}$ and $f \in \CAL{S} '(\R^n)$. Then the short-time Fourier transform $V_{g}f$ of $f$ with respect to a window function $g$ is defined by
\begin{equation*}
V_{g}f(x, \xi)=\int_{\R^n} \overline{g(y-x)}f(y)e^{-iy\cdot \xi}dy,\quad (x, \xi) \in{\mathbb R}^n\times{\mathbb R}^n.
\end{equation*}
We also define the formal adjoint operator $V_{g}^{*}$ of $V_{g}$ by
\begin{equation*}
V_{g}^{*}F(x)=\iint_{{\mathbb R}^{2n}} g(x-y)F(y, \xi)e^{ix\cdot \xi}dy\bar{d}\xi,\quad x\in{\mathbb R}^n, F\in \CAL{S}^{'}(\R^{2n})
\end{equation*}
with $\bar{d}\xi=(2\pi)^{-n}d\xi$. 
\end{Def}
\begin{Prop}[Inverse formula for STFT]\label{Prop1}
Let $g \in \CAL{S}(\mathbb{R}^n)\setminus \{0\}$. Then it follows that
\begin{equation*}
f=\|g\|_{L^2}^{-2} V_{g}^{*}[V_{g}f], \hspace{3mm} f\in \CAL{S}'(\R^n).
\end{equation*}
\end{Prop}
Proposition \ref{Prop1} is proved in \cite[Corollary 11.2.7]{Grochenig}, \cite[Proposition 2.5]{Wahlberg}.
\begin{Prop}[Plancherel's theorem for STFT]\label{Prop2}
Let $g \in \CAL{S}(\mathbb{R}^n)\setminus \{0\}$. Then it follows that
\begin{equation*}
\|V_{g}f\|_{L^{2}(\R^{2n})}=(2\pi)^{n/2} \|g\|_{L^{2}(\R^n)} \|f\|_{L^{2}(\R^{n})},\hspace{3mm} f\in L^{2}(\R^{n}).
\end{equation*}
and
\begin{equation*}
\|V_{g}^{*}F\|_{L^{2}(\R^n)}\leq (2\pi)^{-n/2}\|g\|_{L^{2}(\R^n)} \|F\|_{L^{2}(\R^{2n})},\hspace{3mm} F\in L^{2}(\R^{2n}).
\end{equation*}
\end{Prop}
Proposition \ref{Prop2} is easy to 
be proved by using the Plancherel theorem for $\SCR{F}$ and the Schwartz inequality. Thus we omit the proof here.
\begin{Lem}\label{Lem2}
Let $g \in \CAL{S}(\mathbb{R}^n)\setminus \{0\}, 1\leq p,q \leq \infty$ and $T>0$. We put $g(t)=e^{\frac{1}{2}it\Delta}g$.
If $F(x,\xi)\in L_{x}^{p} L_{\xi}^{q}$, then the estimate holds for all $s,t \in [-T,T]$:
\begin{align*}
\|V_{g(s)}V_{g(t)}^{*}F(x,\xi)\|_{L_{x}^{p}L_{\xi}^{q}}\leq C_{T} \|F(x,\xi)\|_{L_{x}^{p}L_{\xi}^{q}}.
\end{align*}
\end{Lem}
Lemma \ref{Lem2} is proven by an estimate for the solution to the free Schr\"{o}dinger equation and employing the result of Kato-Kobayashi-Ito \cite{KKI1}. The proof is in Appendix.  
\subsection{Lorentz spaces}
We recall the Lorentz space  $L^{p,q}(\mathbb{R}^n)$ for $1\leq p,q\leq \infty$, which is defined by all measurable functions $f$ such that the following quantity is finite:
\begin{align*}
&\|f\|^{*}_{pq}=
\begin{cases}
\left( \frac{q}{p} \int_{0}^{\infty} \bigl[ t^{1/p} f^{*}(t) \bigr]^{q} \frac{dt}{t} \right)^{1/q} \hspace{5mm}\text{if}\hspace{3mm}  1 \leq p < \infty, 1 \leq q < \infty 
\\
\sup_{t>0} t^{1/p} f^{*}(t) \hspace{21mm}\text{if}\hspace{3mm}1 \leq p \leq \infty, q = \infty.
\end{cases}
\end{align*}
Here, $f^{*}$ is the non-increasing rearrangement of $f$ (see \cite{BeLo}).
We have $L^{p,q_{1}} \hookrightarrow L^{p,q_{2}}$ if $q_{1} \leq q_{2}$, and $L^{p,p} = L^{p}$.
Moreover, for  $p\ne1$ the Lorentz space
$L^{p,q}$ becomes a Banach space equipped with the norm $\| \cdot \|_{L^{p,q}}$ which is equivalent to the above quasi-norm $\| \cdot \|^{*}_{pq}$.
\subsection{Wiener amalgam spaces}
We recall the Wiener amalgam spaces introduced by Feichtinger \cite{Fei2}.
We also refer to \cite{Fei1, Fei4, Fournier-Stewart, Heil}. 
First let us define the abstract Wiener amalgam spaces $W(B,C)$ with two Banach  function spaces $B, C$. The Wiener amalgam space $W(B,C)$ is the space of functions which have the local regularity of a function in $B$ and globally display a $C$ decay.
We fix $g\in \CAL{S}(\R^n)\setminus \{0\}$ and put $f_{B}(x)=\|f T_{x}g\|_{B}$. Then we denote
\begin{align*}
&W(B,C):=\{f\in \CAL{S}'(\R^n): \|f\|_{W(B,C)}<\infty\},
\\
&\|f\|_{W(B,C)}:=\|f_{B}\|_{C}.
\end{align*}
As typical examples, it is used as $B=L^q$ or $\SCR{F}L^q$ and $C=L^p$,
The present paper mainly deals with $W(\SCR{F}L^q, L^p)$. Thus we write $W^{p,q}=W(\SCR{F}L^q, L^p)$ for short, whose notation is followed by \cite{Cunanan-Sugimoto}, \cite{Guo-Zhao}.
Then $W^{p,q}$ is defined by the STFT as follows:
\begin{Def}[Wiener amalgam spaces]\label{DefWA}
Let $1\leq p, q\leq \infty$ and $g\in \CAL{S}(\R^n)\setminus \{0\}$. Then we define
\begin{align*}
&W^{p,q}(\R^n)=\left\{f\in\CAL{S}'(\R^n) : \left\| \| V_{g}f(x,\xi)\|_{L_{\xi}^{q}}\right\|_{L_{x}^{p}}<\infty\right\},
\\
&\|f\|_{W^{p,q}(\R^n)}=\left\|  \| V_{g}f(x,\xi)\|_{L_{\xi}^{q}}\right\|_{L_{x}^{p}}.
\end{align*}
\end{Def}
Wiener amalgam spaces are independent of the choice of the window function $g$ (see \cite{Grochenig}).
We collect properties of Wiener amalgam spaces in order to prove Theorem \ref{thm1}.
The following proposition is established by Feichtinger \cite{Fei1, Fei2}, Fournier-Stewart \cite{Fournier-Stewart} and Heil \cite{Heil}.
We also refer to \cite{Cordero-Nicola1, Cordero-Nicola3, Cordero-Nicola2}.
\begin{Prop} \label{Prop4}
Let $B_j,\; C_j, (j = 1,2,3)$ be Banach spaces such that $W(B_j,C_j)$ are well-defined.
Then:
\begin{itemize}
\item[(i)] \textup{(Convolution)}
If $B_1 * B_2 \hookrightarrow B_3$ and $C_1 * C_2 \hookrightarrow C_3$, we have
\begin{equation*}
  W(B_1,C_1) * W(B_2,C_2) \hookrightarrow W(B_3,C_3).
\end{equation*}
In particular, for every $1 \le p,q \le \infty$, we have
\begin{equation*}
\|f * u\|_{W^{p,q}} \leq \|f\|_{W^{1,\infty}} \|u\|_{W^{p,q}}.  
\end{equation*}
\item[(ii)] \textup{(Inclusions)}
If $B_1 \hookrightarrow B_2$ and $C_1 \hookrightarrow C_2$, then
\begin{equation*}
  W(B_1,C_1) \hookrightarrow W(B_2,C_2).
\end{equation*}
In particular, it holds that
\begin{align*}
& L^p \hookrightarrow W^{p',p} \hspace{3mm} \text{if} \hspace{3mm} 1\leq p\leq2,
\\
&W^{p,p'}\hookrightarrow L^p  \hspace{3mm} \text{if} \hspace{3mm} 2\leq p\leq\infty.
\end{align*}
\item[(iii)] \textup{(Complex interpolation)}
For $0 < \theta < 1$, we have
\begin{align*}
  [\,W(B_1,C_1),\, W(B_2,C_2)\,]_{\theta}
  = W\bigl([B_1,B_2]_{\theta},\, [C_1,C_2]_{\theta}\bigr),
\end{align*}
if $C_1$ or $C_2$ has absolutely continuous norm.
In particlar,
for $0<\theta<1$ and $p,q,p_j,q_j\in [1,\infty]$ for $j=1,2$. If $1/p=\theta/p_1+(1-\theta)/p_2$ and $1/q=\theta/q_1+(1-\theta)/q_2$, it holds that
\begin{align*}
\left(W^{p_1,q_1}, W^{p_2,q_2}\right)_{[\theta]}=W^{p,q}.
\end{align*}
\item[(iv)] \textup{(Duality)}
If $B',C'$ are the topological dual spaces of the Banach spaces $B,C$,
respectively, and the space of test functions $C_0^{\infty}$ is dense in both
$B$ and $C$, then
\begin{equation*}
  W(B,C)' = W(B',C'). 
\end{equation*}
In particular, we see that for $1< p,q < \infty$,
\begin{align*}
(W^{p,q})'=W^{p',q'}.
\end{align*}
\end{itemize}
\end{Prop}
\subsection{Classical trajectories}
In this subsection we consider the potential $V$ satisfying Assumption \ref{AssV}. 
We deal with the following Hamilton equation corresponding to (\ref{CP1}):
\begin{align} \label{H}
\begin{cases}
\dot{x}(t)=\xi(t),
\\ \dot{\xi}(t)=-(\nabla_{x}V)(x(t)),
\end{cases}
\end{align}
where the dots are the derivative in time $t$.
The Cauchy–Lipschitz theorem guarantees the global-in-time existence and uniqueness of solutions, and the associated Hamiltonian flow $\Phi(t): \R^{2n}\to \R^{2n}$ enjoys the following properties (see e.g. Fujiwara \cite{Fujiwara1}).
\begin{itemize}
\item $\left(\Phi(t)\right)(x,\xi)$ is the solution to (\ref{H}) with the initial condition $\left(x(0),\xi(0)\right)=(x,\xi)$.
\item $\Phi(0)=$Id.
\item $\Phi(t)\Phi(s)=\Phi(t+s)$\hspace{2mm} for \hspace{2mm}$t,s\in \R$.
\item $\Phi(t)^{-1}=\Phi(-t)$\hspace{2mm} for \hspace{2mm} $t\in \R$.
\end{itemize}
We often write $\left(x(t;x,\xi), \xi(t;x,\xi)\right)=\left(\Phi(t)\right)(x,\xi)$.
Moreover, the conservation of volume in phase space (Liouville's theorem) holds:
\begin{align} \label{Jacobian}
\left|\det \left(\frac{\partial \left(\Phi(t)(x,\xi)\right)}{\partial(x,\xi)}\right) \right|=1
\end{align}
(see e.g. \cite{KKI2}).
\begin{Lem} \label{LemH}
Let $M:=1+n^2 \displaystyle\max_{|\alpha|=2}\displaystyle\sup_{x\in \R^{n}}|\partial_{x}^{\alpha}V(x)|$. We take a constant $T_1>0$ sufficiently small so that $2M^{3/2}e^{MT_1}T_1<1/2$ and $T_1<1/3$ hold. Then we have the inequalities:
\begin{align}
&|x(t;x,\xi)-x(t;z,\eta)|\geq \frac{1}{6}(5|x-z|-3|\xi-\eta|), \label{X(t)}
\\
&|\xi(t;x,\xi)-\xi(t;z,\eta)|\geq \frac{1}{2}(|\xi-\eta|-|x-z|) \label{Xi(t)}
\end{align}
for all $x,z,\xi,\eta\in\R^n$ and $|t|<T_1$.
\end{Lem}
\begin{proof}
We write 
\begin{align} \label{defX}
X(t)=x(t;x,\xi)-x(t;z,\eta), \Xi(t)=\xi(t;x,\xi)-\xi(t;z,\eta), X=x-z \hspace{3mm}\text{and} \hspace{3mm}\Xi=\xi-\eta
\end{align} 
for short. 
Let $E(t):=|X(t)|^2+|\Xi(t)|^2$.
Then (\ref{H}) implies
\begin{align*}
\dot{E}(t)&=2X(t)\cdot \dot{X}(t)+2 \Xi(t)\cdot \dot{\Xi}(t)
\\
&=2X(t)\cdot \Xi(t)+2 \Xi(t)\cdot \{\nabla V(x(t;x,\xi))-\nabla V(x(t;z,\eta))\},
\end{align*}
which together with Taylor's expansion formula of order one and the AM-GM inequality follows that
\begin{align*}
\left|\dot{E}(t)\right|&\leq 2|X(t)\cdot \Xi(t)|+2 n^2 \displaystyle\max_{|\alpha|=2}\displaystyle\sup_{x\in \R^{n}}|\partial_{x}^{\alpha}V(x)| |\Xi(t)| |X(t)|
\\
&\leq 2M E(t).
\end{align*}
By the Gronwall inequality we have
\begin{align*}
 |X(t)|^2+|\Xi(t)|^2 \leq 2M(|X|^2+|\Xi|^2) e^{2M|t|},
\end{align*}
which implies
\begin{align} \label{energyest}
|X(t)|+|\Xi(t)| \leq (4M)^{1/2}(1+|X|+|\Xi|) e^{M|t|}.
\end{align}
In particular, we obtain 
\begin{align} \label{X}
|X(t)|\leq (4M)^{1/2} e^{M|t|} (|X|+|\Xi|).
\end{align}
By the integration of (\ref{H}) and Taylor's theorem we have 
\begin{align}
\Xi(t)&=\Xi-\int_{0}^{t}\left\{\nabla V(x(s;x,\xi))-\nabla V(x(s;z,\eta))\right\}ds \notag
\\
&=\Xi-\int_{0}^{t}\sum_{|\alpha|=1} \left(\int_{0}^{1} (\nabla \partial^{\alpha} V)((x(s;z,\eta)+\theta X(s))\cdot X(s) d\theta\right) ds \label{est1}
\end{align}
and 
\begin{align}
X(t)&=X+\int_{0}^{t}\Xi(s)ds \notag
\\
&=X+t\Xi-\int_{0}^{t}\int_{0}^{s} \sum_{|\alpha|=1} \left(\int_{0}^{1} (\nabla \partial^{\alpha} V)(x(s';z,\eta)+\theta X(s'))\cdot  X(s') d\theta\right) ds' ds. \label{est2}
\end{align}
We find from (\ref{X}), (\ref{est1}) and (\ref{X}), (\ref{est2}) for sufficiently small $T_1$ that
\begin{align}
|\Xi(t)|&\geq |\Xi|-MT_1\cdot (4M)^{1/2} e^{MT_1} (|X|+|\Xi|) \notag
\\
&=\left(1-2M^{3/2} e^{MT_1}T_1\right)|\Xi|-2M^{3/2} e^{MT_1}T_1|X|   \notag
\\
&\geq \frac{1}{2}(|\Xi|-|X|)
\end{align}
and
\begin{align}
|X(t)|&\geq |X|-T_1|\Xi|-MT_{1}^{2} \cdot (4M)^{1/2} e^{MT_1} (|X|+|\Xi|)  \notag
\\
&=\left(1-2M^{3/2} e^{MT_1}T_{1}^{2}\right)|X|-T_1\left(1+2M^{3/2} e^{MT_1}T_1\right) |\Xi|  \notag
\\
&\geq \frac{1}{6}(5|X|-3|\Xi|),
\end{align}
respectively.
\end{proof}
The following Lemma can be obtained by straightforward calculations and using the equation (\ref{H}) and Assumption \ref{AssV} on the potential. We include the proof in Appendix for the reader's convenience.
\begin{Lem} \label{Lemdet}
There exists a small constant $T_2>0$ such that 
\begin{align*}
\frac{1}{2} \leq \left|\det \left(\frac{\partial \left(x(t;x,\xi/t)\right)}{\partial \xi}\right) \right|\leq 2
\end{align*}
for all $(t,x,\xi)\in [-T_2,T_2]\times \R^{2n}$.
\end{Lem}
\section{Key Lemmas}
\begin{Lem} \label{Lem3}
Let $T_1>0$ be as in Lemma \ref{LemH}.
Then there exist some constant $C_{T_1}>0$ such that 
\begin{align} \label{MainInq31}
\|\Phi_{t}V_{g(s)}V_{g(s)}^{*} F(x,\xi)\|_{L_{x}^{\infty}L_{\xi}^{1}} \leq C_{T_1} \|\Phi_{t}F(x,\xi)\|_{L_{x}^{\infty}L_{\xi}^{1}}
\end{align}
for any $s, t\in[-T_1, T_1]$. 
\end{Lem}
\begin{proof}
First we note that for $g\in \CAL{S}, \alpha\in \Z_{+}^{n}$ and $N\in \N$,
\begin{align} \label{freeest}
|\partial_{x}^{\alpha}e^{it\Delta/2} g(x)|\leq C_{N,\alpha,g}\J{x}^{-2N} \J{t}^{2N}.
\end{align}
Indeed, for the $\alpha=0$, we have from integration by parts 
\begin{align*}
&\J{x}^{2N}|e^{it\Delta/2}f(x)|=\left|\int e^{ix\xi} (1-\Delta_{\xi})^{N} [e^{-it|\xi|^2/2} \widehat{g}(\xi)] d\xi\right|
\\
&\leq \sum_{|\beta+\gamma|\leq 2N} \int |\partial_{\xi}^{\beta} e^{-it|\xi|^2/2}| |\partial_{\xi}^{\gamma} \widehat{g}(\xi)| d\xi \leq C_{N,g} \J{t}^{2N}.
\end{align*}
For the $\alpha\neq0$, we can prove (\ref{freeest}) the same argument as the case $\alpha=0$ since $\partial^{\alpha} e^{it\Delta/2}=e^{it\Delta/2} \partial^{\alpha}$.
Integration by parts and (\ref{freeest}) derive
\begin{align*}
&\left|\Phi_{t}V_{g(s)}V_{g(s)}^{*} F(x,\xi)\right|=\left| \int \overline{g(s,y-x(t;x,\xi))} V_{g(s)}^{*} F(y) e^{-iy\xi(t;x,\xi)} dy \right|
\\
&=\left| \iiint \overline{g(s,y-x(t;x,\xi))} g(s,y-z) F(z,\eta) e^{iy(\eta-\xi(t;x,\xi))} dydzd\eta \right|
\\
&\leq C_{N_1} \sum_{|\alpha+\beta|\leq 2 N_1} \iiint |\partial^{\alpha}\overline{g(s,y-x(t;x,\xi))}| |\partial^{\beta}g(s,y-z)| |F(z,\eta)| \J{\eta-\xi(t;x,\xi)}^{-2N_1}dydzd\eta
\\
&\leq C_{N_1, N_2,g, T_1} \iiint \J{y-x(t;x,\xi)}^{-N_2} \J{y-z}^{-N_2-n-1} |F(z,\eta)| \J{\eta-\xi(t;x,\xi)}^{-2N_1}dydzd\eta.
\end{align*}
 Noting that $y-x(t)=\left(z-x(t)\right)+(y-z)$, we have by Peetre's inequality
\begin{align*}
&\left|\Phi_{t}V_{g(s)}V_{g(s)}^{*} F(x,\xi)\right| 
\\
&\leq C_{T_1}  \iiint \J{z-x(t;x,\xi)}^{-N_2} \J{y-z}^{-n-1} |F(z,\eta)| \J{\eta-\xi(t;x,\xi)}^{-2N_1}dydzd\eta
\\
&\leq C_{N, T_1}  \iint \J{z-x(t;x,\xi)}^{-N_2} |F(z,\eta)| \J{\eta-\xi(t;x,\xi)}^{-2N_1}dzd\eta.
\end{align*}
We consider the change of variable : $(z,\eta)\mapsto \Phi_{t}(z,\eta)=\left(x(t;z,\eta), \xi(t;z,\eta)\right)$.
Then (\ref{Jacobian}) yields
\begin{align*}
&\left|\Phi_{t}V_{g(s)}V_{g(s)}^{*} F(x,\xi)\right| 
\\
&\leq C_{N, T_1}  \iint \J{x(t,z,\eta)-x(t;x,\xi)}^{-N_2} |\Phi_{t}F(z,\eta)| \J{\xi(t,z,\eta)-\xi(t;x,\xi)}^{-2N_1}dzd\eta
\\
&= C_{N, T_1}  \iint \J{X(t)}^{-N_2} |\Phi_{t}F(z,\eta)| \J{\Xi(t)}^{-2N_1}dzd\eta.
\end{align*}
Thus we put 
\begin{align} \label{defI}
&I(t,x):=\iiint \J{X(t)}^{-N_2} \J{\Xi(t)}^{-2N_1} |\Phi_{t}F(z,\eta)| dz d\xi d\eta,
\\
&A(x,\xi):=\frac{1}{6}(5|x|-3|\xi|) \hspace{3mm} \text{and} \hspace{3mm}B(x,\xi):=\frac{1}{2}(|\xi|-|x|). \notag
\end{align}
Then the inequality 
\begin{align} \label{inq231}
\|\Phi_{t}V_{g(s)}V_{g(s)}^{*} F(x,\xi)\|_{L_{\xi}^{1}}\leq C_{N,T_1} I(t,x)
\end{align}
holds. In addition we have 
\begin{align*}
&\J{X(t)}^{-1}\leq \J{\max\{A(X,\Xi),0\}}^{-1},
\\
&\J{\Xi(t)}^{-1}\leq \J{\max\{B(X,\Xi),0\}}^{-1}
\end{align*}
for all $t\in[-T_1,T_1]$ by Lemma \ref{LemH}. Hence we consider the following decomposition:
\begin{align*}
I(t,x)=\sum_{i,j\in\{\pm\}}I_{i,j}(t,x)
\end{align*}
where
\begin{align*}
&I_{i,j}(t,x)=\iiint_{D_{i,j}} \J{X(t)}^{-N_2} \J{\Xi(t)}^{-2N_1} |\Phi_{t}F(z,\eta)| dz d\xi d\eta,
\\
&D_{\pm,\pm}=\{(z,\xi,\eta)\in\R^{3n} | \pm A(X,\Xi)\geq0 \hspace{2mm} \& \hspace{2mm} \pm B(X,\Xi)\geq0\}.
\end{align*}
We note that $I_{-,-}(t,x)=0$ because $D_{-,-}$ is the null set in $\R^{3n}$.
Hence we estimate for only $I_{+,+}, I_{+,-}$ and $I_{-,+}$.
We estimate $I_{+,+}$. Translating : $\xi\mapsto \xi-\eta$ we have
\begin{align*} 
&I_{+,+}(t,x)\leq \iint \left( \int \J{A(x-z,\xi-\eta)}^{-N_2} \J{B(x-z,\xi-\eta)}^{-2N_1} d\xi \right) |\Phi_{t}F(z,\eta)| dz d\eta 
\\
&\leq C_{N,N_1} \iiint \J{5|x-z|-3|\xi|}^{-N_2} \J{|\xi|-|x-z|)}^{-2N_1}  |\Phi_{t}F(z,\eta)| dz d\xi d\eta 
\\
&\leq C_{N, N_1}  \|\Phi_{t}F(z,\eta)\|_{L_{z}^{\infty}L_{\eta}^{1}} \iint \J{5|z|-3|\xi|}^{-N_2}  \J{|\xi|-|z|}^{-2N_1}  dz d\xi.  
\end{align*}
We also use the conversion to polar coordinates and obtain
\begin{align} 
&I_{+,+}(t,x)\leq C_{N_1, N_2}  \|\Phi_{t}F(z,\eta)\|_{L_{z}^{\infty}L_{\eta}^{1}} \int_{0}^{\infty} \int_{0}^{\infty} \J{5r-3\rho}^{-N_2} \J{\rho-r}^{-2N_1} r^{n-1} \rho^{n-1} dr d\rho  \notag
\\
&\leq C_{ N, N_1} \|\Phi_{t}F(z,\eta)\|_{L_{z}^{\infty}L_{\eta}^{1}} \int_{0}^{\infty} \int_{-\infty}^{\rho} \J{2\rho-5r'}^{-N_2} \J{r'}^{-2N_1} \J{r'+\rho}^{n-1} \J{\rho}^{n-1} dr' d\rho,  \notag
\end{align}
which implies with the Peetre inequalities
\begin{align}
I_{+,+}(t,x)\leq C_{N_1, N_2} \|\Phi_{t}F(z,\eta)\|_{L_{z}^{\infty}L_{\eta}^{1}} \left(\int \J{\rho}^{-N_2+2(n-1)} d\rho \right) \left(\int \J{r'}^{-2N_1+N_2+n-1} dr' \right). \label{I1}
\end{align}
Next we estimate $I_{+,-}$. By the change of variables : $\xi\mapsto \xi-\eta$ and $z\mapsto x-z$, we have
\begin{align} \label{I2}
&I_{+,-}(t,x)\leq C_{N_1}\iint \left( \int  \chi\{B(x-z,\xi-\eta)\leq0\} \J{A(x-z,\xi-\eta)}^{-N_2} d\xi \right) |\Phi_{t}F(z,\eta)| dz d\eta  \notag
\\
&\leq C_{N_1, N_2}\iint \left( \int  \chi\{B(x-z,\xi)\leq0\} \J{A(x-z,\xi)}^{-N_2} d\xi \right) |\Phi_{t}F(z,\eta)| dz d\eta  \notag
\\
&\leq C_{N_1, N_2}\|\Phi_{t}F(z,\eta)\|_{L_{z}^{\infty}L_{\eta}^{1}} \iint \chi\{B(z,\xi)\leq0\} \J{A(z,\xi)}^{-N_2}  dz d\xi   \notag
\\
&\leq C_{N_1, N_2} \|\Phi_{t}F(z,\eta)\|_{L_{z}^{\infty}L_{\eta}^{1}} \iint \chi\{|\xi|\leq |z|\} \J{5|z|-3|\xi|}^{-N_2}  dz d\xi  \notag
\\
&\leq C_{N_1, N_2} \|\Phi_{t}F(z,\eta)\|_{L_{z}^{\infty}L_{\eta}^{1}} \iint \chi\{|\xi|\leq |z|\} \J{2|z|}^{-N_2}  dz d\xi  \notag
\\
&\leq C_{N_1, N_2} \|\Phi_{t}F(z,\eta)\|_{L_{z}^{\infty}L_{\eta}^{1}} \int \J{z}^{-N_2+n} dz. 
\end{align}
Similarly to (\ref{I2}), we can estimate $I_{-,+}$ as follows:
\begin{align} \label{I3}
&I_{-,+}(t,x)\leq C_{N_2} \iint \left( \int \chi\{A(x-z,\xi-\eta)\leq0\} \J{B(x-z,\xi-\eta)}^{-2N_1} d\xi \right)   |\Phi_{t}F(z,\eta)| dz d\eta  \notag
\\
&\leq C_{N_2} \iint \left( \int \chi\{A(x-z,\xi)\leq0\}  \J{B(x-z,\xi)}^{-2N_1} d\xi \right) |\Phi_{t}F(z,\eta)| dz d\eta  \notag
\\
&\leq C_{N_2} \|\Phi_{t}F(z,\eta)\|_{L_{z}^{\infty}L_{\eta}^{1}} \iint \chi\{A(z,\xi)\leq0\}  \J{B(z,\xi)}^{-2N_1} dz d\xi   \notag
\\
&= C_{N_2} \|\Phi_{t}F(z,\eta)\|_{L_{z}^{\infty}L_{\eta}^{1}} \iint \chi\{5|z|\leq 3|\xi|\} \J{|\xi|-|z|}^{-2N_1} dzd\xi \notag
\\
&\leq C_{N_2} \|\Phi_{t}F(z,\eta)\|_{L_{z}^{\infty}L_{\eta}^{1}} \iint \chi\{5|z|\leq 3|\xi|\} \J{2|\xi|/5}^{-2N_1} dzd\xi \notag
\\
&\leq C_{N_2} \|\Phi_{t}F(z,\eta)\|_{L_{z}^{\infty}L_{\eta}^{1}} \int \J{\xi}^{-2N_1+n} d\xi 
\end{align}
We choose $N_1=3n$ and $N_2=3n$ in order to make the integrals of (\ref{I1})--(\ref{I3}) finite. Therefore by (\ref{inq231})--(\ref{I3}), we have (\ref{MainInq31}).
\end{proof}
\begin{Lem} \label{Lem4}
Let $T_2>0$ be as in Lemma \ref{Lemdet}.
Then we have the estimate
\begin{align} \label{Lem10}
\|\Phi_{t} V_{g(s)}f(x,\xi)\|_{L_{x}^{\infty} L_{\xi}^{1}} \leq C_{T_2} |t|^{-n} \|f\|_{W^{1,\infty}}
\end{align}
for all $s,t\in [-T_2,T_2]$.
\end{Lem}
\begin{proof}
Proposition \ref{Prop1} and the integration by parts show
\begin{align} \label{inq2}
&|\Phi_{t}V_{g(s)}f(x,\xi)|=\left|\int \overline{g(s,y-x(t;x,\xi))} V_{g}^{*} V_{g}f(y) e^{-iy\xi(t;x,\xi)} dy\right| \notag
\\
&= \left|\iiint \overline{g(s,y-x(t;x,\xi))} g(y-z) V_{g}f(z, \eta) e^{iy(\eta-\xi(t;x,\xi))} dy dz d\eta\right|  \notag
\\
&\leq \sum_{|\alpha+\beta|\leq 2n} \iiint \J{\eta-\xi(t;x,\xi)}^{-2n} |\partial^{\alpha}g(s,y-x(t;x,\xi))|   \notag|\partial^{\beta}g(y-z)||V_{g}f(z, \eta)|dy dz d\eta
\\
&\leq C \sum_{|\alpha+\beta|\leq 2n} \iint |\partial^{\alpha}g(s,y-x(t;x,\xi))| |\partial^{\beta}g(y-z)|\|V_{g}f(z, \eta)\|_{L_{\eta}^{\infty}}dy dz.  
\end{align}
By using  (\ref{freeest}), we can further estimate the inequality (\ref{inq2}) as follows:
\begin{align} \label{inq3}
&|\Phi_{t}V_{g(s)}f(x,\xi)|\leq C_{T_2} \sum_{|\beta|\leq 2n} \iint \J{y-x(t;x,\xi)}^{-2n} |\partial^{\beta}g(y-z)|\|V_{g}f(z, \eta)\|_{L_{\eta}^{\infty}}dy dz  
\end{align}
Taking $L_{\xi}^{1}$-norm to both sides of (\ref{inq3}), we can apply Hadamard's inverse mapping theorem (see e.g. \cite{Fujiwara3, Schwartz}) by Lemma \ref{Lemdet} and obtain
\begin{align*}
&\|\Phi_{t}V_{g(s)}f(x,\xi)\|_{L_{\xi}^{1}}\leq C_{T_2} |t|^{-n} \sum_{|\beta|\leq 2n} \iint \left( \int \J{y-x(t;x,\xi/t)}^{-2n}  d\xi \right) |\partial^{\beta}g(y-z)|\|V_{g}f(z, \eta)\|_{L_{\eta}^{\infty}}dy dz 
\\
&\leq C_{T_2} |t|^{-n} \sum_{|\beta|\leq 2n} \iint \left( \int \J{y-\xi}^{-2n}  d\xi \right) |\partial^{\beta}g(y-z)|\|V_{g}f(z, \eta)\|_{L_{\eta}^{\infty}}dy dz 
\\
&\leq C_{T_2} |t|^{-n} \|V_{g}f(z,\eta)\|_{L_{z}^{1} L_{\eta}^{\infty}},
\end{align*}
which implies (\ref{Lem10}).
\end{proof}
\section{Proof of Theorem}
\subsection{Reduction and Transformation}
First we note that the inequality (\ref{MainInq2}) is equivalent to (\ref{MainInq1}) by usual duality argument (see e.g. \cite{Keel-Tao, Cordero-Nicola1, Cordero-Nicola2}).
In addition it suffices to show (\ref{MainInq1}) and (\ref{MainInq3}) for sufficiently small $T>0$.
Indeed, if (\ref{MainInq1}) holds for some $T_0>0$, then  by the facts $U(t+s)=U(t)U(s)$ and $\|U(t)\|_{L^2\to L^2}=$Id, 
\begin{align*}
&\int_{-T}^{T} \|U(t)u_0\|_{W^{p,p'}}^{r/2}dt=\sum_{j=0}^{2N-1}\int_{-T+\frac{T}{N}j}^{-T+\frac{T}{N}(j+1)}\|U(t)u_0\|_{W^{p,p'}}^{r/2}dt
\\
&=\sum_{j=0}^{2N-1}\int_{0}^{\frac{T}{N}}\|U\left(t+T-\frac{T}{N}j\right)u_0\|_{W^{p,p'}}^{r/2}dt
\\
&\leq \sum_{j=0}^{2N-1}\int_{0}^{T_0}\left\|U(t) \left(U(T-\frac{T}{N}j)u_0\right)\right\|_{W^{p,p'}}^{r/2}dt
\\
&\leq \sum_{j=0}^{2N-1}\int_{0}^{T_0}\left\|U\left(T-\frac{T}{N}j\right)u_0\right\|_{L^2}^{r/2}dt=\sum_{j=0}^{N-1}\int_{0}^{T_0}\left\|u_0\right\|_{L^2}^{r/2}dt
\\
&\leq C_{T} \|u_0\|_{L^2}^{r/2}
\end{align*}
with some $N>T/T_0$. Furthermore, if there exists a $T_0>0$ such that (\ref{MainInq3}) is valid, then we obtain for any $T>0$ 
\begin{align*}
&\int_{-T}^{T} \left\|\int_{0}^{t}U(t-s)F(s)ds \right\|_{W^{p,p'}}^{r/2}dt
\\
&=\sum_{j=0}^{2N-1}\int_{-T+\frac{T}{N}j}^{-T+\frac{T}{N}(j+1)} \left\| \int_{0}^{t}U(t-s) F(s) ds \right\|_{W^{p,p'}}^{r/2}dt
\\
&\leq C_r \sum_{j=0}^{2N-1}\int_{-T+\frac{T}{N}j}^{-T+\frac{T}{N}(j+1)} \left\{ \left\| \int_{0}^{-T+\frac{T}{N}j}|U(t-s) F(s)ds \right\|_{W^{p,p'}}^{r/2} \right.
\\
& \hspace{50mm}+ \left. \left\| \int_{-T+\frac{T}{N}j}^{t} U(t-s) F(s)ds \right\|_{W^{p,p'}}^{r/2} \right\}dt,
\end{align*}
which and translating yield
\begin{align} \label{Id41}
&\int_{-T}^{T} \left\|\int_{0}^{t}U(t-s)F(s)ds \right\|_{W^{p,p'}}^{r/2}dt  \notag
\\
&\leq C_r \sum_{j=0}^{2N-1}\int_{0}^{\frac{T}{N}} \left\{ \left\| \int_{0}^{-T+\frac{T}{N}j} U(t-T+\frac{T}{N}j-s) F(s)ds \right\|_{W^{p,p'}}^{r/2} \right.   \notag
\\
&\hspace{30mm}+\left. \left\|\int_{-T+\frac{T}{N}j}^{t-T+\frac{T}{N}j} U(t-T+\frac{T}{N}j-s) F(s)ds \right\|_{W^{p,p'}}^{r/2} \right\}dt   \notag
\\
&= C_r \sum_{j=0}^{2N-1}\int_{0}^{\frac{T}{N}} \left\{ \left\|U(t)U\left(\frac{T}{N}j-T\right)\int_{0}^{-T+\frac{T}{N}j} U(-s) F(s)ds \right\|_{W^{p,p'}}^{r/2} \right.   \notag
\\
&\hspace{30mm}+\left. \left\| \int_{0}^{t} U(t-s) F\left(s+\frac{T}{N}j-T\right)ds \right\|_{W^{p,p'}}^{r/2} \right\}dt
\end{align}
with some $N>T/T_0$.
Applying (\ref{MainInq1}), $\|U(t)\|_{L^2 \to L^2}=1$ and (\ref{MainInq2}) to the first term and (\ref{MainInq3}) for some $T_0>0$ to the second term, we see that (\ref{MainInq3}) for any $T>0$.
 In what follows, we assume $t\in[-T,T]$ with $T<\min\{T_1, T_2\}$, where $T_1, T_2$ are as in Lemmas \ref{LemH}, \ref{Lemdet}, respectively by this reduction.
Furthermore, by a density argument we may assume without loss of generality that $u_0\in \CAL{S}$. Then it holds  that the solution $u\in C^{1}(\R; \CAL{S}(\R^n))$ since $U(t)=e^{-itH}$ is an isomorphism on the Schwartz space (see Fujiwara \cite{Fujiwara1}).
The original equation (\ref{CP1}) can be transformed into
\begin{align}\label{CP1'}
V_{g(t)}[u(t)](x,\xi)&=e^{-i\int_{0}^{t}h(\tau-t;x,\xi)d\tau}V_{g}u_0(x(-t;x,\xi),\xi(-t;x,\xi)) \notag
\\&\hspace{4mm}-i\int_{0}^{t}e^{-i\int_{s}^{t}h(\tau-t;x,\xi)d\tau}\widetilde{R}u\left(s,x(s-t;x,\xi),\xi(s-t;x,\xi)\right)ds,
\end{align}
where $h(t;x,\xi)=\frac{1}{2}|\xi(t;x,\xi)|^2+V(x(t;x,\xi))-\nabla_{x}V(x(t;x,\xi))\cdot x(t;x,\xi)$
and
\begin{align}\label{Ru}
\widetilde{R}u(t,x,\xi)&=\frac{1}{2}\displaystyle\sum_{|\alpha|=2}\int_{\R^n}\int_{0}^{1}\left(\partial_{x}^{\alpha}V\right)(x+\theta(y-x)) (1-\theta)d\theta \notag
\\&\hspace{30mm}\times (y-x)^{\alpha}\overline{g(t, y-x)}u(t,y)e^{-iy\xi}dy
\end{align}
by the method of Kato-Kobayashi-Ito \cite{KKI2}. We give the proof of this transformation for the reader's convenience in Appendix.
 Taking $V_{g(t)}^{*}[\cdot]$ to the both sides of (\ref{CP1'}), we have by Proposition \ref{Prop1},
\begin{align*} 
u(t)= V_{g(t)}^{*} M(t,0) \Phi_{-t} V_{g} u_0 -i\int_{0}^{t} V_{g(t)}^{*} M(t,s) \Phi_{s-t} \widetilde{R}(s) u(s) ds, 
\end{align*}
where $M(t,s)=e^{-i\int_{0}^{t}h(\tau;t,x,\xi)d\tau}$ is the multiplier and $(\Phi_{t}F)(x,\xi)=F\left(\Phi(t)(x,\xi)\right)$ is the change of variables by the Hamilton flow of (\ref{H}). We note that $M(t,s)^{*}=M(s,t)$ and $\Phi_{t}^{*}=\Phi_{-t}$ since $\Phi(t)$ has the group structure and the conservation of volume (\ref{Jacobian}).
We define
\begin{align}
& U_0(t)=V_{g(t)}^{*} M(t,0) \Phi_{-t} V_{g}  \label{U_0},
\\
&R(t,s)=V_{g(t)}^{*} M(t,s) \Phi_{s-t} \widetilde{R}(s) \label{R'}
\end{align}
and obtain
\begin{align} \label{eq1}
U(t)u_0=U_0(t)u_0-i\int_{0}^{t} R(t,s) [U(s)u_0] ds.
\end{align}
\subsection{Estimate of $U_0$}
We have by Lemma \ref{Lem2},
\begin{align} 
&\|U_0(t)U_0(s)^{*}f\|_{W^{\infty,1}}=\|V_{g} \left(V_{g(t)}^{*} M(t,0)\Phi_{-t}V_{g}\right) \left(V_{g}^{*} \Phi_{s}M(0,s) V_{g(s)}\right) f(x,\xi)\|_{L_{x}^{\infty} L_{\xi}^{1}} \notag
\\
&\leq C_{T} \|M(t,0) \Phi_{-t}V_{g} V_{g}^{*} \Phi_{s}M(0,s) V_{g(s)} f(x,\xi)\|_{L_{x}^{\infty} L_{\xi}^{1}} \notag
\\
&= C_{T} \|\Phi_{-t}V_{g} V_{g}^{*} \Phi_{s}M(0,s) V_{g(s)} f(x,\xi)\|_{L_{x}^{\infty} L_{\xi}^{1}}, \notag
\end{align}
which gives with Lemma \ref{Lem3}
\begin{align} \label{inq5}
&\|U_0(t)U_0(s)^{*}f\|_{W^{\infty,1}} \leq C_{T} \| \Phi_{s-t}M(0,s) V_{g(s)} f(x,\xi)\|_{L_{x}^{\infty} L_{\xi}^{1}}  \notag
\\
&= C_{T} \| \Phi_{s-t}V_{g(s)} f(x,\xi)\|_{L_{x}^{\infty} L_{\xi}^{1}}. 
\end{align}
Combining (\ref{inq5}) and Lemma \ref{Lem4}, we obtain
\begin{align} \label{inq6}
\|U_0(t)U_0(s)^{*}f\|_{W^{\infty,1}} \leq C_{T} |s-t|^{-n} \|f\|_{W^{1,\infty}}.
\end{align}
In addition we have by Lemma \ref{Prop2} and (\ref{Jacobian})
\begin{align} \label{inq7}
\|U_0(t)f\|_{W^{2,2}}=\|U_0(t)f\|_{L^2}\leq C_{T} \|f\|_{L^2}=C_{T} \|f\|_{W^{2,2}}.
\end{align}
By Proposition \ref{Prop4} and the same argument as Keel-Tao \cite{Keel-Tao} and Cordero-Nicola \cite{Cordero-Nicola1, Cordero-Nicola2}, we obtain
\begin{align} \label{inq100}
\|U_0(t)f\|_{L^{r}([-T,T];W^{p,p'})}\leq C_{p,T} \|f\|_{L^2}
\end{align}
with $2\leq p \leq \infty,\hspace{2mm} 4< r \leq \infty,\hspace{2mm}\frac{n}{p}+\frac{2}{r}=\frac{n}{2}$
 and 
\begin{align} \label{inq101}
\|U_0(t)f\|_{L^{2}\left([-T,T]; W(\SCR{F}L^{p',2}, L^p)\right)} \leq C_{p,T} \|f\|_{L^2}
\end{align}
with $p=\frac{2n}{n-1}, \hspace{2mm} n>1$.
\subsection{Estimate of $R$}
First we deal with $\widetilde{R}(s)$. By a direct computation, we have
\begin{align*}
&(\widetilde{R}(s)^{*}F)(x)
\\
&=\frac{1}{2}\sum_{|\alpha|=2} \iint (x-y)^{\alpha} g(s,x-y) \left( \int_{0}^{1} (\partial^{\alpha}V)(y+\theta(x-y)) d\theta \right) F(y,\xi) e^{ix\xi} dy d\xi,
\\
&\equiv \frac{1}{2} \sum_{|\alpha|=2} \iint (x-y)^{\alpha} g(s,x-y) V^{\alpha}(x,y) F(y,\xi) e^{ix\xi} dy d\xi
\end{align*}
with $V^{\alpha}\in C_{b}^{\infty}(\R^{2n})$ for $|\alpha|=2$ and 
\begin{align*}
&\widetilde{R}(s)\widetilde{R}(s)^{*}F(x,\xi)= \frac{1}{4} \sum_{|\alpha|=|\beta|=2}\iiint  (y-x)^{\alpha} \overline{g(s,y-x)} (y-z)^{\beta} g(s,y-z) 
\\
& \hspace{50mm} \times V^{\alpha}(y,x) V^{\beta}(y,z) F(z,\eta) e^{iy(\eta-\xi)} dy dz d\eta.
\end{align*}
From integration by parts, the inequality (\ref{freeest}) and Assumption \ref{AssV}, we obtain
\begin{align*}
&|\Phi_{t} \widetilde{R}(s)\widetilde{R}(s)^{*}F(x,\xi)|
\\
&=\frac{1}{4}\left|\sum_{|\alpha|=|\beta|=2}\iiint  (y-x(t;x,\xi))^{\alpha} \overline{g(s,y-x(t;x,\xi))} (y-z)^{\beta} g(s,y-z) V^{\alpha}(y,x(t;x,\xi)) \right.
\\
& \hspace{10mm} \times \left.  V^{\beta}(y,z) F(z,\eta) \J{\eta-\xi(t;x,\xi)}^{-2N_1} (1-\Delta_{y})^{N_1} e^{iy(\eta-\xi(t;x,\xi))} dy dz d\eta \right|
\\
&\leq C_{T} \iiint \J{\eta-\xi(t;x,\xi)}^{-2N_1} \J{y-x(t;x,\xi)}^{-N_2} \J{y-z}^{-N_3} |F(z,\eta)| dy dz d\eta
\end{align*}
for $N_1, N_2>3n$ and $N_3>N_2+n$. By virtue of Peetre's inequality we have
\begin{align*}
&|\Phi_{t} \widetilde{R}(s)\widetilde{R}(s)^{*}F(x,\xi)|
\\
&\leq C_{T} \iiint \J{\eta-\xi(t;x,\xi)}^{-2N_1} \J{z-x(t;x,\xi)}^{-N_2} \J{y-z}^{-N_3+N_2} |F(z,\eta)| dy dz d\eta
\\
&\leq C_{T} \iint \J{\eta-\xi(t;x,\xi)}^{-2N_1} \J{z-x(t;x,\xi)}^{-N_2} |F(z,\eta)|dz d\eta,
\end{align*}
which shows with Liouville's theorem (\ref{Jacobian})
\begin{align*}
&|\Phi_{t} \widetilde{R}(s)\widetilde{R}(s)^{*}F(x,\xi)| \leq C_{T} \iint \J{\Xi(t;x,\xi)}^{-2N_1} \J{X(t;x,\xi)}^{-N_2} |\Phi_{t}F(z,\eta)|dz d\eta
\end{align*}
and 
\begin{align*}
&\|\Phi_{t} \widetilde{R}(s)\widetilde{R}(s)^{*}F(x,\xi)\|_{L_{\xi}^{1}} \leq C_{T} I(t,x),
\end{align*}
where $X(t;x,\xi), \Xi(t;x,\xi)$ and $I(t,x)$ are denoted by (\ref{defX}) and (\ref{defI}).
Hence, we have by the same argument as Lemma \ref{Lem3},
\begin{align} \label{inq8}
\|\Phi_{t} \widetilde{R}(s)\widetilde{R}(s)^{*}F(x,\xi)\|_{L_{x}^{\infty}L_{\xi}^{1}} \leq C_{T} \|\Phi_{t}F(x,\xi)\|_{L_{x}^{\infty}L_{\xi}^{1}}
\end{align}
for all $t_1, t_2, s\in [-T,T]$.
It follows from Lemma \ref{Lem2} and (\ref{inq8}) that 
\begin{align*}
&\|R(t_1,s)R(t_2,s)^{*} f\|_{W^{\infty,1}}
\\
&=\left\| V_{g} \left(V_{g(t_1)}^{*} M(t_1,s) \Phi_{s-t_1} \widetilde{R}(s)\right)\left( \widetilde{R}(s)^{*} \Phi_{t_2-s} M(s,t_2) V_{g(t_2)} \right)f\right\|_{L_{x}^{\infty} L_{\xi}^{1}}
\\
&\leq C_{T} \left\| \Phi_{s-t_1} \Phi_{t_2-s} M(s,t_2) V_{g(t_2)} f(x,\xi)\right\|_{L_{x}^{\infty} L_{\xi}^{1}}
\\
&= C_{T} \left\| \Phi_{t_2-t_1} V_{g(t_2)} f(x,\xi)\right\|_{L_{x}^{\infty} L_{\xi}^{1}}
\end{align*}
which implies with Lemma \ref{Lem4}
\begin{align} \label{inq9}
\|R(t_1,s)R(t_2,s)^{*} f\|_{W^{\infty,1}}\leq C_{T} |t_1-t_2|^{-n} \|f\|_{W^{1,\infty}} 
\end{align}
for all $t_1, t_2, s\in [-T,T]$.
Moreover, we have
\begin{align} \label{inq10}
\|R(t,s)f\|_{W^{2,2}} \leq C_{T} \|f\|_{W^{2,2}}
\end{align}
for all $s,t\in[-T,T]$.
Indeed, we have from Plancherel theorem
\begin{align*}
\left\|\widetilde{R}(s)f(x,\xi)\right\|_{L_{x,\xi}^2}&\leq \frac{1}{2}\displaystyle\sum_{|\alpha|=2}\left\|\int_{0}^{1}\left(\partial_{x}^{\alpha}V\right)(x+\theta(y-x)) (1-\theta)d\theta \right. \notag
\\&\left. \hspace{30mm}\times (y-x)^{\alpha}\overline{g(s, y-x)}f(y) \right\|_{L_{x,y}^2}
\\
&\leq C_{T} \|f\|_{L^2}, 
\end{align*}
which together with (\ref{R'}) yields (\ref{inq10}). 
By Proposition \ref{Prop4} and a similar argument to Keel-Tao \cite{Keel-Tao} and Cordero-Nicola \cite{Cordero-Nicola1, Cordero-Nicola2} with (\ref{inq9}) and (\ref{inq10}),  we have for all $s\in [-T,T]$,
\begin{align} \label{inq11}
\left\| R(\cdot,s)g(s)\right\|_{L^{r}([-T,T];W^{p,p'})}\leq C_{p,T} \|g(s)\|_{L^2},
\end{align}
with $2\leq p \leq \infty,\hspace{2mm} 4< r \leq \infty,\hspace{2mm}\frac{n}{p}+\frac{2}{r}=\frac{n}{2}$
 and 
\begin{align} \label{inq102}
\|R(\cdot,s)g(s)\|_{L^{2}\left([-T,T]; W(\SCR{F}L^{p',2}, L^p)\right)} \leq C_{p,T} \|g(s)\|_{L^2}
\end{align}
with $p=\frac{2n}{n-1}, \hspace{2mm} n>1$.
\subsection{Proof of homogeneous Strichartz estimate}
We write 
\begin{align*}
X(p,r):=
\begin{cases}
&L^r([-T,T];W^{p,p'}) \hspace{2mm} \text{if}\hspace{2mm} 2\leq p \leq \infty,\hspace{2mm} 4< r \leq \infty,\hspace{2mm}\frac{n}{p}+\frac{2}{r}=\frac{n}{2}
\\
&L^{2} \left([-T,T]; W(\SCR{F}L^{p',2}, L^p)\right) \hspace{2mm} \text{if}\hspace{2mm} p=\frac{2n}{n-1}, \hspace{2mm} n>1
\end{cases}
\end{align*}
and $X'(p,r)$ is the topological dual space of $X(p,r)$.
Taking $X(p,r)$-norm to the both sides of the equality (\ref{eq1}), we have by (\ref{inq100}), (\ref{inq101}), (\ref{inq11}) and (\ref{inq102}) 
\begin{align*}
\|u\|_{X(p,r)}&\leq \|U_0(\cdot)u_0\|_{X(p,r)}+ \int_{0}^{T} \left\| \widetilde{R}(\cdot,s)u(s) \right\|_{X(p,r)} ds
\\
&\leq C_{T} \|u_0\|_{L^2}+ C_{T} \int_{0}^{T} \|u(s)\|_{L^2} ds
\\
&\leq C_{T} \|u_0\|_{L^2},
\end{align*}
where we used the conservation law: $\|u(s)\|_{L^2}=\|u_0\|_{L^2}$ and which proves (\ref{MainInq1}) and (\ref{MainInq4}). By virtue of reduction in Subsection 4.1, we also attain (\ref{MainInq2}) and (\ref{MainInq5}).
\subsection{Proof of retarded Strichartz estimate}
In this subsection we prove the inequality (\ref{MainInq3}) and its endpoint case.
 In order to show this, we prepare the following Lemma.
\begin{Lem} \label{LemR}
We have 
\begin{align} \label{id90}
\left(i \partial_{t}-H\right) U_0(t)f=\CAL{R}(t)f,
\end{align}
where $\CAL{R}(t)$ satisfies
\begin{align*}
&\|\CAL{R}(t)f\|_{L^2}\leq C_T \|f\|_{L^2}
\\
&\|\CAL{R}(t) \CAL{R}(s)^{*}f\|_{W^{\infty,1}}\leq C_T \|f\|_{W^{1,\infty}}
\end{align*}
for all $t\in[-T,T]$.
\end{Lem}
We postpone this proof to Appendix.
Let us set
\begin{align*}
\Gamma_0 F(t):=\int_{0}^{t} U_0(t-s)F(s)ds, \hspace{2mm} \Gamma_{\CAL{R}}F(t):=\int_{0}^{t} \CAL{R}(t-s) F(s) ds.
\end{align*}
By the same argument as Keel-Tao \cite{Keel-Tao} and Cordero-Nicola \cite{Cordero-Nicola1, Cordero-Nicola2} we obtain
\begin{align}
&\left\| \Gamma_0 F \right\|_{X(p,r)} \leq C_{T} \|F\|_{X'(\tilde{p},\tilde{r})}, \label{inq451}
\\
&\left\| \Gamma_{\CAL{R}} F \right\|_{X(p,r)} \leq C_{T} \|F\|_{X'(\tilde{p},\tilde{r})} \label{inq452}
\end{align}
for any admissible pairs $(p,r), (\tilde{p}, \tilde{r})$ since $U_0$ has (\ref{inq6}), (\ref{inq7}) and Lemma \ref{LemR}. 
We note that $H=-\frac{1}{2}\Delta+V$ and obtain by Lemma \ref{LemR}
\begin{align*}
&\frac{d}{dt} (e^{itH}U_0(t)f)
\\
&=i e^{itH}\left(-\frac{1}{2}\Delta+V\right) U_0(t)f+ie^{itH}\left\{ \left(\frac{1}{2}\Delta-V\right) U_0(t)f-\CAL{R}(t)f \right\}
\\
&=-i e^{itH} \CAL{R}(t)f.
\end{align*}
We integrate from $0$ to $t$ this equality and obtain 
\begin{align*}
e^{itH} U_0(t) f=f-i \int_{0}^{t} e^{isH} \CAL{R}(s) fds.
\end{align*}
Thus we have
\begin{align*} 
 U_0(t) f= e^{-itH} f-i \int_{0}^{t} e^{-i(t-s)H} \CAL{R}(s) fds,
\end{align*}
which implies with the Fubini theorem
\begin{align} \label{Id110}
\int_{0}^{t} e^{-i(t-s)H} F(s) ds&=\int_{0}^{t} U_0(t-s) F(s) ds +i \int_{0}^{t} \int_{0}^{t-s} e^{-i(t-s-\tau)H} \CAL{R}(\tau) F(s) d\tau ds \notag
\\
&=\Gamma_0 F(t) +i \int_{0}^{t} \int_{s}^{t} e^{-i(t-\tau)H} \CAL{R}(\tau-s) F(s) d\tau ds \notag
\\
&=\Gamma_0 F(t) +i \int_{0}^{t}  e^{-i(t-\tau)H} \left[\int_{0}^{\tau} \CAL{R}(\tau-s) F(s) ds\right] d\tau \notag
\\
&=\Gamma_0 F(t) +i \int_{0}^{t}  e^{-itH} \left[e^{i\tau H} \Gamma_{\CAL{R}} F(\tau)\right] d\tau.
\end{align}
We take $X(p,r)$-norm to both sides of the identity (\ref{Id110}) and obtain by (\ref{MainInq1}) and (\ref{inq451})
\begin{align*}
&\left\|\int_{0}^{t} e^{-i(t-s)H} F(s) ds\right\|_{X(p,r)}
\\
&\leq \|\Gamma_0 F\|_{X(p,r)}+\int_{0}^{T} \left\| e^{-itH} \left[e^{i\tau H} \Gamma_{\CAL{R}} F(\tau)\right] \right\|_{X(p,r)} d\tau
\\
&\leq  C_T \|F\|_{X'(\tilde{p},\tilde{r})}+ C_T \int_{0}^{T} \left\|e^{i\tau H} \Gamma_{\CAL{R}} F(\tau) \right\|_{L^2} d\tau
\end{align*}
 Moreover, applying the conservation law and (\ref{inq452}) to this inequality we have
\begin{align*}
\left\|\int_{0}^{t} e^{-i(t-s)H} F(s) ds\right\|_{X(p,r)}&\leq C_T \| F\|_{X'(\tilde{p},\tilde{r})}+\int_{0}^{T} \left\| \Gamma_{\CAL{R}} F \right\|_{X(2,\infty)} d\tau
\\
&\leq C_T \| F\|_{X'(\tilde{p},\tilde{r})},
\end{align*}
which completes the proof.
\section*{Appendix}
\appendix
\subsection*{Proof of Lemma \ref{Lemdet}}
By the integration of (\ref{H}) we have
\begin{align*}
&\xi(t;x,\xi)=\xi-\int_{0}^{t} (\nabla V)(x(s;x,\xi))ds,
\\
&x(t;x,\xi)=x+\int_{0}^{t}\xi(s)ds=x+t\xi-\int_{0}^{t}(t-s) (\nabla V)(x(s;x,\xi))ds
\end{align*} 
for all $(t,x,\xi)\in \R^{1+2n}$.
Setting $\widetilde{x}(t;x,\xi):=x(t;x,\xi/t)$, we have
\begin{align*}
\widetilde{x}(t;x,\xi)=x+\xi-\int_{0}^{t}(t-s)(\nabla V)(\widetilde{x}(s;x,\xi))ds,
\end{align*}
from which we obtain
\begin{align} \label{eq18}
\partial_{\xi_k}\widetilde{x_j}(t;x,\xi)=\delta_{j,k}-\int_{0}^{t}(t-s)\sum_{l=1}^{n}(\partial_{x_l}\partial_{x_j}V)(\widetilde{x}(s;x,\xi)) \partial_{\xi_k} \widetilde{x_l}(s;x,\xi) ds.
\end{align}
Let $M':=1+\displaystyle\max_{|\alpha|=2}\displaystyle\sup_{x\in \R^{n}}|\partial_{x}^{\alpha}V(x)|$. Then  we have by (\ref{eq18})
\begin{align*}
&\|\partial_{\xi_k}\widetilde{x}(\cdot\hspace{1mm} ;x,\xi)\|_{L^{\infty}[-T_2,T_2]}\leq 1+nM' T_{2}^{2} \|\partial_{\xi_k}\widetilde{x}(\cdot\hspace{1mm};x,\xi)\|_{L^{\infty}[-T_2,T_2]},
\end{align*}
that is,
\begin{align} \label{xi1}
\|\partial_{\xi_k}\widetilde{x}(\cdot\hspace{1mm};x,\xi)\|_{L^{\infty}[-T_2,T_2]}\leq \frac{1}{1-nM' T_{2}^{2}}, 
\end{align}
where $\|\partial_{\xi_k}\widetilde{x}(\cdot\hspace{1mm} ;x,\xi)\|_{L^{\infty}[-T_2,T_2]}=\max_{1\leq j\leq n}\|\partial_{\xi_k}\widetilde{x_j}(\cdot\hspace{1mm};x,\xi)\|_{L^{\infty}[-T_2,T_2]}$.
Combining (\ref{eq18}) and (\ref{xi1}) imply that for any $1\leq j, k\leq n$, 
\begin{align} \label{xi2}
|\partial_{\xi_k}\widetilde{x_j}(t,x,\xi)| \geq 1-\frac{n M' T_{2}^{2}}{1-nM'T_{2}^{2}}. 
\end{align}
Hence choosing small $T_{2}>0$, we have the desired assertion by (\ref{xi1}) and (\ref{xi2}).
\subsection*{Proof of the transformation (\ref{CP1}) into (\ref{CP1'})}
First we note that the solution for $u\in C^{1}(\R; \CAL{S}(\R^n))$ when $u_0\in \CAL{S}(\R^n)$.
We take $g(t)=e^{\frac{1}{2}it\Delta}g$ with $g\in{\CAL S}({\R}^n)\setminus \{0\}$ and $\|g\|_{L^2}=1$ as the window function of STFT. 
We have by the Leibniz rule
\begin{align} \label{A1}
V_{g(t)}[i\partial_{t}u(t)]=i\partial_{t}V_{g(t)}[u(t)]+V_{i\partial_{t}g(t)}[u(t)]
\end{align} 
and obtain from integration by parts
\begin{align} \label{A2}
V_{g(t)}[\Delta u(t)]&=\int u(t,y) \Delta_{y}[\overline{g(t,y-x)} e^{-iy\xi}]dy \notag
\\
&=\int u(t,y) \{\Delta_{y}\overline{g(t,y-x)}-2i\xi\cdot \nabla_{y}\overline{g(t,y-x)}-\overline{g(t,y-x)}|\xi|^2\} e^{-iy\xi} dy  \notag
\\
&=\int u(t,y) \{\Delta_{y}\overline{g(t,y-x)}+2i\xi\cdot \nabla_{x}\overline{g(t,y-x)}-\overline{g(t,y-x)}|\xi|^2\} e^{-iy\xi} dy   \notag
\\
&=V_{\Delta g(t)}[u(t)]+2i\xi\cdot \nabla_{x} V_{g(t)}[u(t)]-|\xi|^2 V_{g(t)}[u(t)].
\end{align}
The Taylor's expansion
\begin{align*}
V(y)=V(x)+(y-x)\cdot \nabla V(x)+\frac{1}{2}\sum_{|\alpha|=2}\int_{0}^{1}\left(\partial_{x}^{\alpha}V\right)(x+\theta(y-x)) (1-\theta)d\theta (y-x)^{\alpha}
\end{align*}
follows that
\begin{align} \label{A3}
V_{g(t)}[Vu(t)]&=V(x)V_{g(t)}u(t)+ \int \overline{g(t,y-x)}u(t,y)\nabla V(x) \cdot ye^{-iy\xi}dy \notag
\\
&\hspace{5mm}-x\cdot \nabla V(x) V_{g(t)}u(t)+Ru(t)  \notag
\\
&=V(x)V_{g(t)}u(t)+i \nabla V(x)\cdot \nabla_{\xi} V_{g(t)}u(t)-x\cdot \nabla V(x) V_{g(t)}u(t)+Ru(t),
\end{align}
where $Ru(t)$ is defined by (\ref{Ru}).
Taking the STFT $V_{g(t)}[\cdot]$ to the both sides of (\ref{CP1}), we have
\begin{align} \label{A4}
&i(\partial_{t}+\xi\cdot \nabla_{x}-\nabla V(x)\cdot \nabla_{\xi})V_{g(t)}[u(t)](x,\xi)+V_{\left(i\partial_{t}+\frac{1}{2}\Delta\right) g(t)}[u(t)] \notag
\\
&=\left(\frac{1}{2}|\xi|^2+V(x)-\nabla V(x) \cdot x \right)V_{g(t)}[u(t)](x,\xi)+Ru(t,x,\xi)
\end{align}
by (\ref{A1}), (\ref{A2}) and (\ref{A3}). We note that $V_{\left(i\partial_{t}+\frac{1}{2}\Delta\right) g(t)}[u(t)]=0$ since $g(t)$ solves the free Schr\"{o}dinger equation.
Let $v(s)=V_{g(s)}[u(s)](x(s-t;x,\xi),\xi(s-t,x,\xi)$, where $\left(x(t;x,\xi), \xi(t;x,\xi)\right)=\left(\Phi(t)\right)(x,\xi)$ is the classical trajectory for (\ref{H}). Then, (\ref{H}) and (\ref{A4}) provide
\begin{align*}
\frac{d}{ds}v(s)&=\left(\partial_{s}V_{g(s)}[u(s)]\right)(x(s-t;x,\xi),\xi(s-t;x,\xi))
\\
&\hspace{5mm}+\dot{x}(s-t;x,\xi)\cdot \nabla_{x}V_{g(s)}[u(s)](x(s-t;x,\xi),\xi(s-t;x,\xi))
\\
&\hspace{5mm}+\dot{\xi}(s-t;x,\xi) \cdot \nabla_{\xi}V_{g(s)}[u(s)](x(s-t;x,\xi),\xi(s-t;x,\xi))
\\
&=\left(\partial_{s}V_{g(s)}[u(s)]\right)(x(s-t;x,\xi),\xi(s-t;x,\xi))
\\
&\hspace{5mm}+\xi(s-t;x,\xi)\cdot \nabla_{x}V_{g(s)}[u(s)](x(s-t;x,\xi),\xi(s-t;x,\xi))
\\
&\hspace{5mm}-\nabla V(x(s-t;x,\xi))\cdot \nabla_{\xi}V_{g(s)}[u(s)](x(s-t;x,\xi),\xi(s-t;x,\xi))
\\
&=-ih(s-t,x,\xi) v(s)-iRu(s,x(s-t;x,\xi),\xi(s-t;x,\xi)), 
\end{align*}
which implies
\begin{align*}
v(t)=e^{-i\int_{0}^{t}h(\tau-t,x,\xi)d\tau}v(0)-i\int_{0}^{t}e^{-i\int_{s}^{t}h(\tau-t,x,\xi)} Ru(s,x(s-t;x,\xi),\xi(s-t;x,\xi))ds,
\end{align*}
that is, we attain (\ref{CP1'}).
In particular, for $V\equiv0$ we obtain
\begin{align} \label{A5}
V_{g(t)}[e^{it\Delta/2}u_0](x,\xi)=e^{-it|\xi|^2/2}V_{g}u_0(x-t\xi,\xi),
\end{align}
which has been shown by Kato-Kobayashi-Ito \cite{KKI1}.
\subsection*{Proof of Lemma \ref{Lem2}}
From the definition of STFT and Young's inequality, we have
\begin{align*}
\|V_{g(s)} V_{g(t)}^{*}F(x,\xi)\|_{L_{x}^{p}L_{\xi}^{q}}&=\|\left(V_{g(s)}[g(t)]*F\right)(x,\xi)\|_{L_{x}^{p}L_{\xi}^{q}}
\\
&\leq C \||V_{g(s)}[g(t)](x,\xi)\|_{L_{x}^{1}L_{\xi}^{1}} \||F(x,\xi)\|_{L_{x}^{p}L_{\xi}^{q}},
\end{align*}
which implies with (\ref{A5}) that
\begin{align} \label{A6}
\|V_{g(s)} V_{g(t)}^{*}F(x,\xi)\|_{L_{x}^{p}L_{\xi}^{q}} \leq C \||V_{g(s-t)}[g](x,\xi)\|_{L_{x}^{1}L_{\xi}^{1}} \|F(x,\xi)\|_{L_{x}^{p}L_{\xi}^{q}}.
\end{align}
Integration by parts and (\ref{freeest}) show
\begin{align*}
|V_{g(s-t)}[g](x,\xi)|&\leq C \sum_{|\alpha+\beta|\leq 2n} \J{\xi}^{-2n} \int |\partial_{y}^{\alpha}\overline{g(s-t,y-x)}| |\partial_{y}^{\beta}g(y)|dy
\\
&\leq C_{T}  \sum_{|\alpha+\beta|\leq 2n} \J{\xi}^{-2n} \int \J{y-x}^{-n-1} |\partial_{y}^{\beta}g(y)|dy,
\end{align*}
which implies 
\begin{align} \label{A7}
\|V_{g(s-t)}[g](x,\xi)\|_{L_{x}^{1}L_{\xi}^{1}}\leq C_{T}.
\end{align}
It follows from (\ref{A6}) and (\ref{A7}) that the assertion.
\subsection*{Proof of Lemma \ref{LemR}}
We have by (\ref{Jacobian})
\begin{align*}
U_0(t)f(x)&=\iint g(t,y-x) e^{-i\int_{0}^{t}h(\tau-t;y,\eta)d\tau}V_{g}f(x(-t;y,\eta),\xi(-t;y,\eta)) e^{ix\eta} dy d\eta
\\
&=\iint g(t,x(t;y,\eta)-x) e^{-i\int_{0}^{t}h(\tau;y,\eta)d\tau}V_{g}f(y,\eta) e^{ix\cdot\xi(t;y,\eta)} dy d\eta.
\end{align*}
We write $A:=e^{-i\int_{0}^{t}h(\tau;y,\eta)d\tau}V_{g}f(y,\eta) e^{ix\cdot\xi(t;y,\eta)}$ for short.
 We can compute by (\ref{H})
\begin{align} \label{A10}
\partial_{t}U_0(t)f(x)=&\iint \left\{(\partial_{t}g)(t,x(t;y,\eta)-x)+\xi(t;y,\eta)\cdot (\nabla_{x}g)(t,x(t;y,\eta)-x) \right\} A\hspace{2mm} dyd\eta \notag
\\
&-i\iint g(t,x(t;y,\eta)-x) \left(\frac{1}{2}|\xi(t)|^2+V(x(t))-x(t)\cdot \nabla V(x(t)) \right) A \notag \hspace{2mm} dy d\eta
\\
&-i\iint g(t,x(t;y,\eta)-x) \hspace{2mm}x\cdot \nabla V(x(t)) A \hspace{2mm} dy d\eta.
\end{align}
Similarly, we can calculate
\begin{align} \label{A11}
&\frac{1}{2}\Delta U_0(t)f(x)=\iint \left( \frac{1}{2}(\Delta g)(t,x(t)-x)-i\xi(t)\cdot (\nabla g)(t,x(t)-x)-\frac{1}{2}|\xi(t)|^2 \right) A \hspace{2mm}dyd\eta \notag
\\
&=\iint \left( -i(\partial_{t} g)(t,x(t)-x)-i\xi(t)\cdot (\nabla g)(t,x(t)-x)-\frac{1}{2}|\xi(t)|^2 \right) A \hspace{2mm}dyd\eta
\end{align}
since $g(t)=e^{it\Delta/2}g$ solves the free Schr\"{o}dinger equation, 
where $x(t)=x(t;y,\eta)$ and $\xi(t)=\xi(t;y,\eta)$.
Moreover we have by the Taylor theorem
\begin{align} \label{A12}
&V(x)U_0(t)f(x)=\iint g(t,x(t;y,\eta)-x) \left\{V(x(t))+(x-x(t))\cdot \nabla V(x(t))\right\}  A\hspace{2mm} dy d\eta  \notag
\\
&+\frac{1}{2} \iint g(t,x(t;y,\eta)-x) \sum_{|\alpha|=2}\int_{0}^{1}\left(\partial_{x}^{\alpha}V\right)(x(t)+\theta(x-x(t)) (1-\theta)d\theta (x-x(t))^{\alpha} A \hspace{2mm} dy d\eta.
\end{align}
Combining (\ref{A10}), (\ref{A11}) and (\ref{A12}), we attain
\begin{align}
&(i\partial_{t}-H)U_0(t)f(x)=\frac{1}{2} \iint g(t,x(t;y,\eta)-x) 
\\
& \times \sum_{|\alpha|=2} \left(\int_{0}^{1}\left(\partial_{x}^{\alpha}V\right)\left(x(t;y,\eta)+\theta(x-x(t;y,\eta))\right) (1-\theta)d\theta \right) (x-x(t;y,\eta))^{\alpha} A \hspace{2mm} dy d\eta.
\end{align}
Let $\CAL{R}(t)f(x)$ be the right hand side of this equality.
Then it holds that
\begin{align*}
&\CAL{R}(t)=\widetilde{ \CAL{R} } (t) \Phi_{-t} M(t,0) V_{g},
\\
&\widetilde{ \CAL{R} }(t)F(x)=\frac{1}{2}\sum_{|\alpha|=2} \iint \overline{g(t,y-x)} (x-y)^{\alpha} \left(\int_{0}^{1} (\partial^{\alpha}V)(y+\theta(x-y)) d\theta \right) F(y,\eta) e^{ix\eta} dyd\eta.
\end{align*}
Consequently, we attain the desired result by an argument similar to Section 4.3.
\subsection*{Proof of the statement in Remark \ref{Rem}}
It suffices to prove $\|e^{it(\Delta/2-E\cdot x)}f\|_{W^{\infty,1}}=\|e^{it\Delta/2}f\|_{W^{\infty,1}}$ from
 the results in the case without a potential by \cite{Cordero-Nicola2} and the conservation law: $\|e^{it(\Delta/2-E\cdot x)}f\|_{W^{2,2}}=\|f\|_{W^{2,2}}$.
It is known that the explicit formula holds:
\begin{align*}
e^{it(\Delta/2-E\cdot x)}f(x)&=(2\pi it)^{-n/2} \int_{\R^n} e^{i\left( \frac{(x-y)^2}{2t}-\frac{1}{2} E(x+y)t -\frac{1}{24}E^2 t^3\right)} f(y) dy
\end{align*}
(see \cite{AH, Ozawa}).
By a simple computation this is rewritten as
\begin{align*}
e^{it(\Delta/2-E\cdot x)}&=e^{iE^2 t^3/6} e^{-itEx} T_{-Et^2/2} e^{it\Delta/2},
\end{align*}
where $T_{a}$ is the translation: $(Tf)(x)=f(x-a)$.
We have by the definition of STFT  
\begin{align*} 
V_{g}[e^{-itE(\cdot)}T_{a}f](x,\xi)
&=e^{-itEa} e^{-ia\xi} V_{g}f(x-a, \xi+tE), 
\end{align*}
which gives
\begin{align*}
\|e^{it(\Delta/2-E\cdot x)}f\|_{W^{\infty,1}}&=\|V_{g}[e^{-itE(\cdot)}  T_{-Et^2/2} e^{it\Delta/2}f](x,\xi)\|_{L_{x}^{\infty} L_{\xi}^{1}}
\\
&=\|V_{g}[e^{it\Delta/2}f](x+Et^2/2, \xi+tE)\|_{L_{x}^{\infty} L_{\xi}^{1}}
\\
&=\|e^{it\Delta/2}f\|_{W^{\infty,1}}
\end{align*}
and completes the proof.
\subsection*{Acknowledgments}
The author would like to thank Professor Yuusuke Sugiyama for valuable discussions and comments on the proof of endpoint case. 

\address{
Shun Takizawa\\
Department of Mathematics,\\
 Faculty of Science, \\
Tokyo University of Science,\\
Kagurazaka 1-3, Shinjuku-ku, \\
Tokyo 162-8601, Japan}
{1123703@ed.tus.ac.jp}

\end{document}